\documentclass[11pt,a4paper]{article}

\usepackage[a4paper,margin=30mm]{geometry}
\usepackage[T1]{fontenc}
\usepackage{lmodern}
\usepackage{microtype}

\usepackage[reqno]{amsmath}
\usepackage{amssymb,amsthm,mathtools,comment}
\usepackage{aliascnt}

\usepackage{enumitem}
\usepackage[numbers,sort&compress]{natbib}
\usepackage{tikz}
\usetikzlibrary{arrows.meta}
\usepackage[colorlinks=true,linkcolor=blue,citecolor=blue,urlcolor=blue]{hyperref}
\usepackage[noabbrev]{cleveref}
\makeatletter
\renewcommand{\section}{\@startsection{section}{1}\z@
  {2.5ex\@plus 1ex\@minus .2ex}{1.2ex\@plus .2ex}
  {\normalfont\large\bfseries\centering}}
\renewcommand{\subsection}{\@startsection{subsection}{2}\z@
  {2ex\@plus .8ex\@minus .2ex}{.8ex\@plus .2ex}
  {\normalfont\normalsize\bfseries\raggedright}}
\makeatother

\renewenvironment{abstract}{%
  \par\medskip
  {\centering\normalfont\normalsize\bfseries Abstract\par}%
  \vspace{0.6em}%
  \begingroup
  \normalfont\normalsize
  \leftskip=3pc\relax
  \rightskip=3pc\relax
  \parindent=1.5em\relax
  \indent\ignorespaces
}{%
  \par\endgroup\medskip
}

\newtheorem{theorem}{Theorem}[section]

\newaliascnt{lemma}{theorem}
\newtheorem{lemma}[lemma]{Lemma}
\aliascntresetthe{lemma}
\newaliascnt{proposition}{theorem}
\newtheorem{proposition}[proposition]{Proposition}
\aliascntresetthe{proposition}
\newaliascnt{corollary}{theorem}
\newtheorem{corollary}[corollary]{Corollary}
\aliascntresetthe{corollary}
\newaliascnt{claim}{theorem}

\aliascntresetthe{claim}
\newaliascnt{conjecture}{theorem}
\newtheorem{conjecture}[conjecture]{Conjecture}
\aliascntresetthe{conjecture}
\newaliascnt{definition}{theorem}

\aliascntresetthe{definition}
\newaliascnt{fact}{theorem}

\aliascntresetthe{fact}
\theoremstyle{remark}
\newaliascnt{remark}{theorem}

\aliascntresetthe{remark}

\crefname{theorem}{Theorem}{Theorems}
\crefname{lemma}{Lemma}{Lemmas}
\crefname{claim}{Claim}{Claims}
\crefname{section}{Section}{Sections}
\crefname{corollary}{Corollary}{Corollaries}
\crefname{proposition}{Proposition}{Propositions}
\crefname{conjecture}{Conjecture}{Conjectures}
\crefname{fact}{Fact}{Facts}
\crefname{equation}{Equation}{Equations}
\crefname{definition}{Definition}{Definitions}
\crefname{remark}{Remark}{Remarks}

\allowdisplaybreaks

\title{On Win's Conjecture for Edge-Disjoint Perfect Matchings}
\author{
Yufei Chang\thanks{School of Mathematics,
Shandong University, Jinan 250100, China.}
\and
Yangyang Cheng\thanks{Corresponding author.
Email: {\tt yangyang.cheng@uni-passau.de}. Faculty of Computer Science
and Mathematics, University of Passau, Passau,
Germany. Partially supported by the Deutsche
Forschungsgemeinschaft (DFG, German Research
Foundation)-542321564.}
\and
Zhilan Wang\footnotemark[1]
\and
Shuo Wei\footnotemark[1]
}
\date{}

\begin{document}
\maketitle

\begin{abstract}
In 1982, Win conjectured that, for every even integer $n$ and every integer $k$
$1\le k\le n-2$, any graph $G$ of order $n$ satisfying
$d_G(x)+d_G(y)\ge n+k-2$ for every pair of nonadjacent
vertices $x,y$ contains $k$ pairwise edge-disjoint
perfect matchings. In this paper, we prove a stronger
result for all sufficiently large even $n$. More precisely,
we show that under the same degree-sum condition,
$G$ contains $\lfloor k/2\rfloor$ pairwise edge-disjoint
Hamilton cycles, together with an additional perfect
matching when $k$ is odd, all pairwise edge-disjoint.
Consequently, we confirm Win's conjecture for all
sufficiently large even $n$, with a threshold independent
of $k$.


\end{abstract}

\section{Introduction}\label{sec:introduction}

Graph decomposition problems form a central and active area of research in graph theory. A \emph{decomposition} of a graph $G$
is a collection of pairwise edge-disjoint subgraphs whose edge sets partition $E(G)$. If all these subgraphs are isomorphic to a prescribed graph $H$, then the decomposition is called an \emph{$H$-decomposition}. 
The study of graph decompositions has a long history.
A classical result of Walecki \cite{WaleckiSurvey}
shows that the complete graph $K_{2m+1}$ can be decomposed into $m$ Hamilton cycles, while $K_{2m}$ can be decomposed into $m-1$ Hamilton cycles and a perfect matching. 

More recently, Csaba, K\"uhn, Lo, Osthus and
Treglown~\cite{CKLOT} resolved the $1$-factorization
and Hamilton decomposition conjectures for sufficiently large graphs. In particular, they
established the following result.

\begin{theorem}[Csaba et al.~\cite{CKLOT}]\label{Ck}
There exists an integer $n_0$ such that every
$r$-regular graph $G$ on $n\ge n_0$ vertices with
$r\ge\lfloor n/2\rfloor$ admits a decomposition
into $\lfloor r/2\rfloor$ Hamilton cycles and,
if $r$ is odd, one perfect matching.
\end{theorem}

Another important result is the Hamilton decomposition
theorem of K\"uhn and Osthus~\cite{KO} for robust expanders.
For $\nu>0$ and $S\subseteq V(G)$, the
\emph{robust $\nu$-neighborhood} of $S$ is defined by
$\operatorname{RN}_{G}(S)=
\{v\in V(G):|N_G(v)\cap S|\ge\nu n\}$.
A graph $G$ of order $n$ is called a
\emph{robust $(\nu,\tau)$-expander} if
$|\operatorname{RN}_{G}(S)|\ge |S|+\nu n$
for every $S\subseteq V(G)$ with
$\tau n\le |S|\le(1-\tau)n$.
K\"uhn and Osthus established the following result.

\begin{theorem}[K\"uhn and Osthus~\cite{KO}]
\label{thm:decomposition}
For every $\alpha>0$, there exist constants
$\tau>0$ and $\nu_0>0$ such that, for every
$0<\nu\le\nu_0$, there exists an integer $n_0$
with the following property. Every even-regular
robust $(\nu,\tau)$-expander $G$ on $n\ge n_0$
vertices with $\delta(G)\ge\alpha n$
admits a Hamilton decomposition.
\end{theorem}

This theorem is a key ingredient in the proof of
\cref{Ck}. Indeed, for every fixed $\varepsilon>0$,
every graph $G$ with
$\delta(G)\ge(1/2+\varepsilon)n$ is a robust
$(\nu,\tau)$-expander. Consequently,
\cref{thm:decomposition} implies the asymptotic
version of \cref{Ck}.

A natural relaxation of a graph decomposition problem
is to require only that the prescribed subgraphs are
pairwise edge-disjoint, without insisting that they cover all edges of the host graph.
A fundamental problem in this direction is to determine sufficient conditions for
the existence of many pairwise edge-disjoint perfect
matchings.
For a graph $G$, let $\sigma_2(G)$ denote the minimum
of $d_G(x)+d_G(y)$ over all pairs of distinct
nonadjacent vertices $x,y$, with $\sigma_2(K_n)=+\infty$.
Win~\cite{Win} proposed the following conjecture for edge-disjoint perfect matchings.

\begin{conjecture}[Win~\cite{Win}]\label{conj:win}
Let $n$ be even and $1\le k\le n-2$. If an $n$-vertex
graph $G$ satisfies $\sigma_2(G)\ge n+k-2$, then $G$
contains $k$ pairwise edge-disjoint perfect matchings.
\end{conjecture}

The cases $k=1,2$ follow from Ore’s theorem \cite{Ore}, while Win~\cite{Win} established the case
$k=3$. For every fixed even $k$, the Hamilton cycle
packing theorem of Faudree, Rousseau and
Schelp~\cite{FRS} also implies \cref{conj:win} for
sufficiently large $n$. A substantial advance was
made by Lu and Ning~\cite{LuNing}, who confirmed the
conjecture for all $k\ge n/2$ when $n$ is sufficiently
large. We state their result explicitly.

\begin{theorem}[Lu and Ning~\cite{LuNing}]\label{Lu}
There exists an integer $n_0$ such that, whenever
$n\ge n_0$ is even and $n/2\le k\le n-2$,
every $n$-vertex graph $G$ satisfying
$\sigma_2(G)\ge n+k-2$ contains $k$ pairwise
edge-disjoint perfect matchings.
\end{theorem}


In this paper, we confirm \cref{conj:win} for all
sufficiently large even $n$ independent of $k$. In fact, we establish the following stronger result, which guarantees edge-disjoint Hamilton cycles together with, when
necessary, one additional perfect matching.

\begin{theorem}\label{thm:main}
There exists an integer $n_0$ such that, whenever
$n\ge n_0$ is even, $1\le k\le n-2$, and $G$ is an
$n$-vertex simple graph satisfying
$\sigma_2(G)\ge n+k-2$, the following assertions hold.
If $k$ is even, then $G$ contains $k/2$ pairwise
edge-disjoint Hamilton cycles. If $k$ is odd, then
$G$ contains $(k-1)/2$ pairwise edge-disjoint
Hamilton cycles together with a perfect matching
edge-disjoint from all of them.
\end{theorem}

Since every Hamilton cycle in an even-order graph
can be decomposed into two perfect matchings,
\cref{thm:main} immediately yields the following
confirmation of \cref{conj:win} for large $n$.

\begin{corollary}\label{cor:win}
There exists an integer $n_0$ such that, for every
even $n\ge n_0$ and every integer $1\le k\le n-2$,
each $n$-vertex graph $G$ with
$\sigma_2(G)\ge n+k-2$ contains $k$ pairwise
edge-disjoint perfect matchings.
\end{corollary}

\section{Preliminaries}\label{sec:preliminaries}
\subsection{Notation}
All graphs are finite, simple and undirected. For a vertex $v$ of
$G$, let $N_G(v)$ and $d_G(v)$ denote its neighborhood and degree,
respectively. 
We write $N_G[v]=N_G(v)\cup\{v\}$ for the
closed neighborhood of $v$. The minimum and maximum degrees are
$\delta(G)=\min_{v\in V(G)}d_G(v)$ and
$\Delta(G)=\max_{v\in V(G)}d_G(v)$. A \emph{vertex cut}, or simply a
\emph{cut}, of $G$ is a set $S\subseteq V(G)$ such that $G-S$ is
disconnected or has at most one vertex, the vertex connectivity
$\kappa(G)$ is the minimum size of such a set. 
For any $A,B\subseteq V(G)$, let
$e_G(A,B)=\sum_{v\in A}|N_G(v)\cap B|$.
In particular, edges with both endpoints in
$A\cap B$ are counted twice.
A \emph{$k$-factor} is a spanning
$k$-regular subgraph.

\subsection{Useful tools}

We begin by recalling several results that will be used in our proofs.

\begin{lemma}[Egawa~\cite{Egawa}]\label{thm:egawa}
Let $r\ge2$ and $n\ge44(r-1)$.  If an $n$-vertex graph $G$ satisfies
$\sigma_2(G)\ge n$ and $\delta(G)\ge4r-2$, then $G$ contains $r$
pairwise edge-disjoint Hamilton cycles.
\end{lemma}
\begin{lemma}[Li~\cite{Li}]\label{thm:li}
Let $r\ge1$ and let $G$ be a $2r$-connected graph on  $n$ vertex.
If every pair of vertices $x,y$ at distance two satisfies
$\max\{d_G(x),d_G(y)\}\ge n/2+2r-2$, then $G$ contains $r$
pairwise edge-disjoint Hamilton cycles.
\end{lemma}

\begin{lemma}[Bondy and Chv\'atal~\cite{BondyChvatal}]
\label{lem:closure}
Let $G$ be an $n$-vertex graph, and let $x,y$ be
nonadjacent vertices of $G$. If
$d_G(x)+d_G(y)\ge n$, then $G$ is Hamiltonian
if and only if $G+xy$ is Hamiltonian.
\end{lemma}

\begin{lemma}[Lu and Ning~\cite{LuNing}]
\label{ln}
Let $n/2\le k\le n-1$ with $kn$ even.
If $G$ is an $n$-vertex graph satisfying
$\sigma_2(G)\ge n+k-2$, then $G$
contains a $k$-factor.
\end{lemma}

For a function $f:V(G)\to\mathbb N$,
an $f$-factor of $G$ is a spanning subgraph $F$
such that $d_F(v)=f(v)$ for every $v\in V(G)$.
We also use Tutte's $f$-factor theorem \cite{Tutte}. For a graph $G$ and $f:V(G)\to\mathbb N$, write
$f(X)=\sum_{v\in X}f(v)$. 

\begin{lemma}[Tutte~\cite{Tutte}]\label{thm:tutte}
A graph $G$ has an $f$-factor if and only if, for all disjoint sets $S,T\subseteq V(G)$, $f(S)+\sum_{v\in T}\bigl(d_{G-S}(v)-f(v)\bigr)\ge q(S,T)$, where $q(S,T)$ count the components $C$ of $G-S-T$ for which
$f(V(C))+e_G(V(C),T)$ is odd.
\end{lemma}

We first establish a connectivity property.
The first assertion of the following proposition
was also proved by Lai and Yan~\cite{LaiYan}.

\begin{proposition}\label{lem:connectivity}
Let \(1\le k\le n-2\), and let \(G\) be an \(n\)-vertex graph satisfying \(\sigma_2(G)\ge n+k-2\). Then \(\kappa(G)\ge k\). Moreover, if \(n>2k\) and \(\kappa(G)=k\), there is an edge
$ab\notin E(G)$ such that $\kappa(G+ab)\ge k+1$.
\end{proposition}
\begin{proof}
Suppose that $S$ is a vertex cut of $G$.
Let $A,B$ be two distinct components of $G-S$.
For any $x\in A$ and $y\in B$, we have
\begin{equation}\label{1}
    n+k-2\le d_G(x)+d_G(y)\le |A|+|B|+2|S|-2
\le n+|S|-2.
\end{equation}
Thus $|S|\ge k$, giving $\kappa(G)\ge k$.

If $\kappa(G)=k$ and $|S|=k$.
Then equality in \eqref{1} holds throughout, so $G-S$ consists of exactly two cliques $G[A]$ and $G[B]$,
and every vertex of $S$ is adjacent to every
vertex of $A\cup B$.

Assume $n>2k$. Choose $a\in A$ and $b\in B$,
and let $G_1=G+ab$. For any $T\subseteq V(G)$
with $|T|\le k$, if $T=S$, then $G_1-T$ is
connected by the edge $ab$. Otherwise,
$S\setminus T\ne\varnothing$ and
$(A\cup B)\setminus T\ne\varnothing$, since
$|A\cup B|=n-k>k$. As all edges between
these two nonempty sets are present,
$G_1-T$ is connected. Hence
$\kappa(G_1)\ge k+1$.
\end{proof}

\begin{lemma}[Chernoff bound \cite{J}]\label{lem:chernoff}
Let $X$ be a sum of independent Bernoulli random
variables with $\mu=\mathbb E[X]$. Then the following holds:

\textnormal{($i$)} For every $\delta>0$, $\mathbb P(X\ge(1+\delta)\mu)
\le\exp\left(-\frac{\delta^2\mu}{2+\delta}\right)$

\textnormal{($ii$)} For every $0<\delta<1$, $\mathbb P(X\le(1-\delta)\mu)
\le\exp\left(-\frac{\delta^2\mu}{2}\right)$.
\end{lemma}

\section{Proof of main Theorem}

\subsection{Small $k$ case.}

The next lemma deals with the range of small $k$.
\begin{lemma}\label{prop:small}
Let $n$ be even, $1\le k\le n-2$, and $n\ge22k$. Every $n$-vertex
graph $G$ with $\sigma_2(G)\ge n+k-2$ satisfies the conclusion of
\cref{thm:main}.
\end{lemma}
\begin{proof}
For $k=1$, add a universal vertex to $G$. The resulting $(n+1)$-vertex graph meets Ore's condition and has a Hamilton cycle; deleting the new vertex yields a Hamilton path in $G$, which contains a perfect matching since $n$ is even. For $k=2$, Ore's theorem supplies a Hamilton cycle directly.

Now suppose $k\ge3$ and let $r=\lceil k/2\rceil$.
We first show that if $G$ is $2r$-connected then it contains $r$ pairwise edge-disjoint Hamilton cycles.

    Starting from $G:=G_1$, we can add an edge $xy\notin E(G_i)$ whenever $d_{G_i}(x)+d_{G_i}(y)\ge n+4r-4$. This process terminates with a  sequence $G_1\subseteq G_2\subseteq \dots\subseteq G_t$ such that every pair of nonadjacent vertices $x,y$ in $G_t$ satisfies
$d_{G_t}(x)+d_{G_t}(y)<n+4r-4$.

We claim that adding edge preserves the existence of $r$ edge-disjoint Hamilton cycles. Indeed, suppose there are $r$ edge-disjoint Hamilton cycles $C_1,\dots,C_r$ in $G_i+xy$, and at most one of its cycles uses $xy$, say $C_r$. Let $H=G_i-\cup_{i<r}E(C_i)$. Then $H+xy$ is Hamiltonian and $d_{H}(x)+d_{H}(y)\ge n+4r-4-4(r-1)\ge n$, by \cref{lem:closure}, $H$ has a Hamilton cycle $C_r'$, combining with $C_1,\dots,C_{r-1}$, there are $r$ edge-disjoint Hamilton cycles in $G_i$.

 It remains to find $r$ edge-disjoint Hamilton cycles in $G_t$. If $\delta(G_t)\ge4r-2$, then \cref{thm:egawa} applies, since $\sigma_{2}(G_t)\ge n+k-2$. Otherwise, choose $x\in V(G_t)$ with $d_{G_t}(x)=d\le 4r-3$, and let $B=V(G_t)\setminus N_{G_t}[x]$. Each $b\in B$ satisfies
$d_{G_t}(b)\ge n+k-2-d>n-d$. Now let $v$ be any vertex with
$d_{G_t}(v)\ge n/2$. If $v$ is nonadjacent
to some $b\in B\setminus\{v\}$, then $d_{G_t}(v)+d_{G_t}(b)\ge \frac n2+n-d\ge n+4r-4,$
where the last inequality follows from
$n\ge44(r-1)$ and $d\le4r-3$.
This contradicts the construction of $G_t$.
Thus $v$ is adjacent to every vertex of
$B\setminus\{v\}$, and hence $d_{G_t}(v)\ge |B|-1=n-d-2\ge \frac n2+2r-2.$
Since $\sigma_2(G_t)\ge n$, every pair of
vertices at distance two has at least one
endpoint of degree at least $n/2$.
Consequently, $\max\{d_{G_t}(u),d_{G_t}(v)\}
\ge \frac n2+2r-2$ for every pair $u,v$ at distance two. As $G_t$ is $2r$-connected, \cref{thm:li}
yields $r$ edge-disjoint Hamilton cycles
in $G_t$, and hence also in $G$.

By \cref{lem:connectivity},
$\kappa(G)\ge k$. If $k=2r$, then $G$ is
$2r$-connected, and the preceding argument
gives $r$ edge-disjoint Hamilton cycles in $G$.
Assume that $k=2r-1$, then $\kappa(G)\ge 2r-1$.
If $\kappa(G)=2r-1$, by \cref{lem:connectivity},
there exist an edge $e\notin E(G)$ such that
$G'=G+e$ is $2r$-connected, otherwise, set
$G'=G$. Applying the preceding argument again,
we obtain $r$ edge-disjoint Hamilton cycles in
$G'$. At most one of these cycles uses $e$,
if added. Deleting $e$ from that cycle if
necessary (or any edge otherwise) yields a
Hamilton path in $G$ edge-disjoint from the
remaining $r-1$ Hamilton cycles. Since $n$
is even, this path contains a perfect matching,
as required.s
\end{proof}

\subsection{Extension to regular factors}

\begin{proposition}\label{ab}
Let $G$ be an $n$-vertex graph satisfying
$\sigma_2(G)\ge n+k-2$. For any subsets
$X,Y\subseteq V(G)$ with $|X|+|Y|>0$, we have
\begin{equation*}
e_G(X,Y)\ge
\frac{|X||Y|(|X|+|Y|-n+k)}{|X|+|Y|}
-|X\cap Y|.
\end{equation*}
\end{proposition}

\begin{proof}
The inequality is trivial if either $X$ or $Y$
is empty. Thus we may assume $|X|,|Y|>0$.
Let $H$ be the bipartite graph with vertex
classes given by disjoint copies of $X$ and $Y$,
in which there is an edge between $x\in X$
and $y\in Y$ if and only if $x\ne y$ and
$xy\notin E(G)$. Consequently, for every
$xy\in E(H)$,
$d_H(x)+d_H(y)\le
d_{\overline G}(x)+d_{\overline G}(y)
\le n-k$. Then, we have
\begin{equation*}
\frac{|E(H)|^2}{|X|}+\frac{|E(H)|^2}{|Y|}
\le \sum_{v\in V(H)}d_H(v)^2
=\sum_{xy\in E(H)}(d_H(x)+d_H(y))
\le(n-k)|E(H)|.
\end{equation*}
Hence $|E(H)|\le(n-k)|X||Y|/(|X|+|Y|)$.
Since
$|X||Y|=e_G(X,Y)+|E(H)|+|X\cap Y|$,
the claim follows.
\end{proof}

\begin{lemma}\label{lem:extension}
For every $0<\alpha\le1/2$, there exist $\varepsilon>0$ and $n_0=n_0(\alpha)$ such that the following holds. Let $n\ge n_0$, $\alpha n\le k\le n/2$ with $kn$ even, and let $G$ be an
$n$-vertex graph satisfying $\sigma_2(G)\ge n+k-2$.
Then every $H\subseteq G$ with $\Delta(H)\le\varepsilon n$ is contained in a $k$-factor of $G$.
\end{lemma}

\begin{proof}
Choose constants
\begin{equation}\label{1}
    0<1/n_0\ll\varepsilon\ll c\ll\eta\ll\alpha\le 1/2.
\end{equation}
Let $H\subseteq G$ with $\Delta(H)\le\varepsilon n$,
and put $G'=G-E(H)$. We claim that
$\kappa(G')\ge k-2\varepsilon n$.
By \cref{lem:connectivity},
$\delta(G')\ge\delta(G)-\Delta(H)\ge k-\varepsilon n$.
Suppose that $G'-W$ is disconnected for some
$W\subseteq V(G)$ with $|W|<k-2\varepsilon n$.
For every component $C$ of $G'-W$ and every
$v\in C$, we have
$d_{G'}(v)\le |C|-1+|W|$,
which implies $|C|>\varepsilon n$.
Let $C_1,C_2$ be two distinct components of $G'-W$.
Since $e_{G'}(C_1,C_2)=0$,
every edge of $G$ between $C_1$ and $C_2$
must belong to $H$. If all edges between
$C_1$ and $C_2$ were present in $G$, then
$d_H(x)\ge |C_2|>\Delta(H)$ for any $x\in C_1$,
a contradiction. Thus there exist $x\in C_1$
and $y\in C_2$ such that $xy\notin E(G)$.
 Consequently, $d_G(x)+d_G(y)
\le |C_1|+|C_2|+2|W|+2\varepsilon n-2
\le n+|W|+2\varepsilon n-2
<n+k-2,$
contrary to the Ore condition. Hence
$\kappa(G')\ge k-2\varepsilon n$.

Define $f(v)=k-d_H(v)$ for every $v\in V(G)$.
Then $f(v)\ge0$ and
$f(V(G))=kn-2e(H)$ is even.
Let $S,T$ be disjoint subset of $V(G)$, we have
\begin{equation}\label{left}
    \begin{aligned}
    &f(S)+\sum_{v\in T}\bigl(d_{G'-S}(v)-f(v)\bigr)\\
    &=k(|S|-|T|)-d_H(S)+d_{H}(T)+e_{G'}(T,V(G)\setminus S)\\
    &=k(|S|-|T|)-d_H(S)+e_{G}(T,V(G)\setminus S)+e_{H}(S,T)\\
  &\ge k(|S|-|T|)+e_{G}(T,V(G)\setminus S)-\varepsilon n|S|\\
\end{aligned}
\end{equation}
Suppose first that $|S|+|T|\le\eta n$.
Since $\kappa(G')\ge k-2\varepsilon n>\eta n$, $G'-S-T$ is connected, so $q(S,T)\le1$. By $\delta(G)\ge k$, we have $e_{G}(T,V(G)\setminus S)\ge |T|(k-|S|)$. Thus, the left-hand side of \eqref{left} is at least $|S|(k-|T|-\varepsilon n)\ge 0$. Since $f(V(G))$ is even and 
\[
f(S)+\sum_{v\in T}\bigl(d_{G'-S}(v)-f(v)\bigr)
\equiv q(S,T)\pmod2
\]
implies that $f(S)+\sum_{v\in T}\bigl(d_{G'-S}(v)-f(v)\bigr)\ge q(S,T)$.
Now suppose that $|S|+|T|>\eta n$. Let $h=|S|-|T|$. 
If $h>k/2$, then $k(|S|-|T|)+e_{G}(T,V(G)\setminus S) \ge k^2/2\ge cn^2$.
Otherwise, $h\le k/2$. Applying
\cref{ab} and using
$|T|(n-|S|)=|S|(n-|T|)-nh$, we obtain $k(|S|-|T|)+e_{G}(T,V(G)\setminus S)$ is greater than
\[
\frac{(n-k)h^2+(k-h)|S|(n-|T|)}{n-h}-|T|\ge \frac{nh^2/2+k|S|^2/2}{2n}-n\ge \frac{\alpha}{20}(|S|+|T|)^2-n\ge cn^2-n.
\]
By \eqref{1} and $\alpha n\le k\le n/2$.
Thus, in all cases,
$k(|S|-|T|)+e_G(T,V(G)\setminus S)\ge cn^2-n$.
Consequently, by \eqref{left}, $f(S)+\sum_{v\in T}\bigl(d_{G'-S}(v)-f(v)\bigr)
\ge(c-\varepsilon)n^2-n
>n\ge q(S,T).$
Thus \cref{thm:tutte} holds for all
disjoint sets $S,T\subseteq V(G)$ and $G'$ has an $f$-factor. Thus, $G$ has $k$-factor containing all the edge of $H$.
\end{proof}

\subsection{Construction of robust expanders}

\begin{lemma}\label{lem:backbone}
For every $0<\alpha\le1/2$ and
$0<\varepsilon,\tau<1/2$, there exist
$\nu>0$ and $n_0$ such that the following holds.
 Let $n\ge n_0$, $\alpha n\le k\le n/2$, and let $G$ be an
$n$-vertex graph satisfying
$\sigma_2(G)\ge n+k-2$.
Then, for every matching $M$ in $G$,
there exists a spanning robust $(\nu,\tau)$-expander
$R\subseteq G-E(M)$ with $\Delta(R)\le\varepsilon n/2$.
\end{lemma}

\begin{proof}
Choose constants such that
$1/n_0\ll\nu\ll \alpha,\tau,\varepsilon\le 1/2$.
Select each edge of $G-E(M)$ independently with
probability $\varepsilon/8$, and let $R$ be the resulting
spanning graph. For every $v\in V(G)$, $\mathbb E(d_R(v))\le \varepsilon n/8$.
By \cref{lem:chernoff}, 
\begin{equation}\label{mun}
    \mathbb P\left(d_R(v)\ge\varepsilon n/2\right)\le \exp\left(-9\varepsilon n/40\right)\le\exp\left(-\nu n\right).
\end{equation}

Let $A,B\subseteq V(G)$ satisfy $|A|,|B|\ge\tau n/2$ and $|A|+|B|\ge(1-\alpha/2)n$. Since $|A||B|/(|A|+|B|)\ge\tau n/4$ and
$|A|+|B|-n+k\ge\alpha n/2$, by \cref{ab}, and noting that removing a matching deletes at most $n$ edges, we obtain
\[
e_{G-E(M)}(A,B)\ge\frac{|A||B|(|A|+|B|-n+k)}{|A|+|B|}-2n\ge\frac{\alpha\tau}{16}n^2.
\]
at least $\alpha\tau n^2/32\ge 2\nu n^2$ distinct edges are counted by $e_{G-E(M)}(A,B)$.
Let $X$ be the number of distinct edges counted by
$e_{G-E(M)}(A,B)$ that are selected into $R$.
Since each edge is counted at most twice,
by the choice of $\nu$, we have $\mathbb E[X]= e_{G-E(M)}(A,B)\cdot\varepsilon/16\ge8\nu n^2$.
By \cref{lem:chernoff} (ii), we have 
\begin{equation}\label{n}
  \mathbb P\left(e_R(A,B)<2\nu n^2\right)\le\mathbb P\left(X<2\nu n^2\right)\le\mathbb P\left(X<\mathbb E(X)/2\right)\le\exp\left(
-\nu n^2\right). 
\end{equation}
Since there are at most $4^n$ ordered pairs
$(A,B)$, a union bound together with \eqref{mun} and \eqref{n} gives $n\exp(-\nu n)+4^n\exp\left(-\nu n^2\right)<1$.
Hence there exists a spanning subgraph
$R\subseteq G-E(M)$ satisfying
$\Delta(R)\le\varepsilon n/2$ and
$e_R(A,B)\ge2\nu n^2$ for every pair
$A,B$ satisfying $A,B\subseteq V(G)$ with $|A|,|B|\ge\tau n/2$ and $|A|+|B|\ge(1-\alpha/2)n$.

Suppose that $R$ is not a robust
$(\nu,\tau)$-expander. Then there exists
$X\subseteq V(G)$ with
$\tau n\le|X|\le(1-\tau)n$ and
$|\operatorname{RN}_{R}(X)|<|X|+\nu n$.
Let $Y=V(G)\setminus\operatorname{RN}_{R}(X)$.
Since $\nu\ll\tau,\alpha$, we have
$|X|,|Y|\ge\tau n/2$ and $|X|+|Y|>(1-\alpha/2)n$.
Hence $2\nu n^2\le e_R(X,Y)<\nu n|Y|\le\nu n^2,$
a contradiction. Thus $R$ is a robust
$(\nu,\tau)$-expander.
\end{proof}

\subsection{Proof of \cref{thm:main}}\label{sec:main-proof}

Let $\tau,\nu_0>0$ be given by \cref{thm:decomposition}. Choose constants
$1/n_0\ll\nu\ll\varepsilon\ll\tau\ll\alpha\le 1/22$ such that $\nu\le \nu_0$.
If $k\le\alpha n$, then
$n\ge22k$, and the conclusion follows
from \cref{prop:small}.
Next suppose that $\alpha n<k<n/2$.
If $k$ is odd, let $M$ be a perfect matching
of $G$, which exists by Ore's theorem.
Otherwise, set $M=\varnothing$.
By \cref{lem:backbone}, there exists a spanning
robust $(\nu,\tau)$-expander
$R\subseteq G-E(M)$ with
$\Delta(R)\le\varepsilon n/2$.
Since $\Delta(R\cup M)\le\varepsilon n$, \cref{lem:extension}
gives a $k$-factor $F$ of $G$ containing
$R\cup M$.
Now $F-E(M)$ is an even-regular graph
containing $R$, and hence is a robust
$(\nu,\tau)$-expander. Moreover, its degree is
$2\lfloor k/2\rfloor\ge\alpha n/2$.
Thus, by \cref{thm:decomposition},
$F-E(M)$ has a Hamilton decomposition into
$\lfloor k/2\rfloor$ pairwise edge-disjoint
Hamilton cycles. Together with $M$ when
$k$ is odd, these give the desired conclusion.

Finally, suppose that $k\ge n/2$.
By \cref{ln}, $G$ contains a $k$-factor $F$. Since $k\ge n/2$, \cref{Ck} gives a
decomposition of $F$ into
$\lfloor k/2\rfloor$ Hamilton cycles and,
when $k$ is odd, a perfect matching.
This completes the proof.

\vspace{0.4cm}
\noindent \textbf{Declaration of AI.}
The authors used GPT-6 Pro as a research assistance tool during the
development of this work. In particular, the authors provided the AI with the paper \cite{CKLOT}. With the assistance of GPT-6 Pro, an initial proof of Win's conjecture was obtained. Building on this proof, the authors subsequently established the stronger
result presented as \cref{thm:main}.

The authors have carefully verified all mathematical arguments
and proofs, and have independently revised and reorganized the
manuscript to improve its mathematical precision, logical
coherence, and readability. All results and proofs have been
reviewed and approved by the authors, who take full responsibility
for the correctness, originality, and final content of this paper.

\end{document}